\documentclass[leqno,12pt]{article}
\usepackage{amsmath,amsfonts,amssymb,amsthm}
\usepackage{amscd}
\usepackage{graphicx}
\usepackage{hyperref}%
\newtheorem{thm}{Theorem}[subsection]

\newtheorem{thmsec}{Theorem}[section]
\newtheorem{crl}{Corollary}

\newtheorem{lem}{Lemma}

\newtheorem{prop}{Proposition}

\theoremstyle{remark}
\newtheorem{rem}{Remark}[subsection]

\theoremstyle{definition}

\newtheorem{defin}{Definition}

 \newcommand{\set}[2]{\{{#1}:{#2}\}}
\newcommand{\rarrowsim}{\smash{\mathop{\,\longrightarrow\,}\limits
  ^{\lower1.5pt\hbox{$\scriptstyle\sim$}}}}
\newcommand{\rarrowsimshort}{\smash{\mathop{\,\rightarrow\,}\limits
  ^{\lower1.5pt\hbox{$\scriptstyle\sim$}}}}

\newcommand{\contra}[1]{{{#1}^\vee}}

\newcommand{\OV}{{\mathcal{O}(D,\mathcal{V})}}    
\newcommand{\barOV}{{\overline{\mathcal{O}}(D,\overline{\mathcal{V}})}}
\newcommand{\itm}[1]{\newline\noindent{\rm{#1}}\enspace}
\newcommand{\End}{\operatorname{End}}

\newcommand{\Aut}{\operatorname{{Aut}}}

\providecommand{\bysame}{\makebox[3em]{\hrulefill}\thinspace}
\newenvironment{texteqn}
{\begin{equation} 
\addtolength{\abovedisplayskip}{-1ex}
\addtolength{\abovedisplayshortskip}{-1ex}
\addtolength{\belowdisplayskip}{-1ex}
\addtolength{\belowdisplayshortskip}{-1ex}
\begin{minipage}[t]{0.87\linewidth}}
{\end{minipage} \end{equation} \ignorespacesafterend}
\newenvironment{texteqn*}
{\begin{equation*} \hphantom{(4.1.2)}
\addtolength{\abovedisplayskip}{-1ex}
\addtolength{\abovedisplayshortskip}{-1ex}
\addtolength{\belowdisplayskip}{-1ex}
\addtolength{\belowdisplayshortskip}{-1ex}
\begin{minipage}[t]{0.87\linewidth}}
{\end{minipage} \end{equation*} \ignorespacesafterend}
{\begin{equation} \hfill\begin{minipage}[t]{0.88\linewidth}}
{\end{minipage} \end{equation} \ignorespacesafterend}

\title{%
Propagation of multiplicity-freeness property\\
for holomorphic vector bundles\footnote{This work was partially supported by Grant-in-Aid for Scientific Research 
18340037, 22340026, 
Japan Society for the Promotion of Science.}
\\
\textit{\large Dedicated to Joseph Wolf for his seventy-fifth birthday}
}
\author{Toshiyuki KOBAYASHI\\[\smallskipamount]
Graduate School of Mathematical Sciences, and IPMU,\\
The University of Tokyo\\
\normalsize
\textit{E-mail address}: \texttt{toshi@ms.u-tokyo.ac.jp}
}
\date{}
\begin{document}
\maketitle

\begin{abstract}
We give a complete proof of a propagation theorem of multiplicity-free property 
from fibers to spaces of global sections for 
holomorphic vector bundles.
The propagation theorem is formalised in three ways,
aiming for producing various multiplicity-free theorems
 in representation theory
 for both finite and infinite dimensional cases in a systematic and
synthetic manner.

               The key geometric condition in our theorem is an 
orbit-preserving
               anti-holomorphic diffeomorphism on the base space,
               which brings us to the concept of visible actions on complex 
manifolds.
\end{abstract}

\noindent
\textit{%
Mathematics Subject Classifications\/} (2010): 
Primary:
          22E46, 
\\
Secondary
32M10, 
32M05, 
46E22. 

\medskip

\noindent
\textit{%
Key words\/}: 
multiplicity-free representation, reproducing kernel, unitary 
representation,
  homogeneous space, holomorphic bundle, visible action

\setcounter{tocdepth}{1}
\tableofcontents

\section{Introduction}

\numberwithin{equation}{section}

In representation theory,
 unitarity is an important concept, in particular, when we
 apply the classic philosophy --- analysis and synthesis, namely,
an attempt to understand things built up from the smallest ones.
This is embodied by
a theorem of Mautner and Teleman stating that
any unitary representation $\pi$ of a locally compact group $G$ can be
decomposed into the direct integral of irreducible unitary
representations: 
\begin{equation}\label{eqn:MTdecomp}
\pi \simeq \int_{\widehat{G}}^\oplus m_\pi(\tau) \tau d\mu (\tau),
\end{equation}
where $\widehat{G}$ denotes the set of equivalence classes of irreducible
unitary representations (`\textit{smallest objects}'),
$\mu$ is a measure on $\widehat{G}$, 
and $m: \widehat{G} \to \mathbb{N} \cup \{\infty\}$
is a measurable function 
that stands for `multiplicity'. 

To find elements of $\widehat{G}$,
basic results are unitarizability theorems
established by Mackey \cite{xmackey}
 for $L^2$-induced representations in the 1950s and
by Vogan \cite{xvu} and Wallach \cite{xwal} for cohomologically 
induced representations in the 1980s.
These results may be thought of as a 
\textit{propagation theory of unitarity}
 from fibers to spaces of sections
(more generally, stalks to cohomologies).

Multiplicity-freeness is another important concept in representation
theory that generalizes irreducibility. 
For a unitary representation $\pi$ of $G$,
we say that $\pi$ is \textit{multiplicity-free}  
if the ring of continuous $G$-intertwining endomorphisms is
commutative. 
This condition implies that $m$ is not
greater than $1$ almost everywhere with respect to the measure $\mu$
in the direct integral \eqref{eqn:MTdecomp}.

Multiplicity-free representations are a special class of
representations, 
for which one could expect a simple and detailed study,
 and by which one could expect effective applications 
of representation theory.  
Multiplicity-free representations arise in broad range of mathematics
in connection with expansions
(Taylor series, Fourier expansion, spherical harmonics,
the Gelfand--Tsetlin basis, \dots)
and the classical identities (the Capelli identity, various explicit formulae 
for special functions, \dots),
though we may not be aware of even the fact that the representation is
there. 

The aim of this paper is to prove a \textit{propagation theorem of
multiplicity-freeness} from fibers to spaces of sections for
holomorphic vector bundles.

To state our main result, 
let $H$ be a Lie group,
and $\mathcal{V} \to D$ an $H$-equivariant holomorphic vector
bundle. 
We naturally have a representation of $H$ on the space
$\mathcal{O}(D, \mathcal{V})$  of global holomorphic sections.
Then, the first form of our multiplicity-free theorem is stated
briefly as follows
(see Theorem~\ref{thm:2.2} for details):

\begin{thmsec}
\label{thm:main}
Any unitary representation of $H$ which is realized in 
$\mathcal{O}(D,\mathcal{V})$ 
is multiplicity-free if the $H$-equivariant bundle
$\mathcal{V}\to D$ satisfies the following 
three conditions:
\begin{texteqn}
{\upshape(}Fiber{\upshape)} \ 
For every $x \in D$,
the isotropy representation of $H_x$ on the fiber $\mathcal{V}_x$ 
is multiplicity-free.
\end{texteqn}
\begin{texteqn}
\label{eqn:1.2}
{\upshape(}Base space{\upshape)} \ 
There exists an anti-holomorphic bundle endomorphism  $\sigma$,
which preserves every $H$-orbit on the base space $D$.
\end{texteqn}
\begin{texteqn}
\label{eqn:1.3}
{\upshape(}Compatibility{\upshape)} \ 
See \eqref{thmassA}.
\end{texteqn}
\end{thmsec}

The compatibility condition \eqref{eqn:1.3} 
is less important because it is often automatically fulfilled by a natural 
choice of $\sigma$
(see Remark~\ref{rem:Weylinv} for an example of $\sigma$;
see also \cite[Appendix]{xkmfkorea}).
Thus, 
 the geometric condition \eqref{eqn:1.2} on the base space $D$ is crucial
 for our propagation of multiplicity-free property from fibers
$\mathcal{V}_x$ to the space
$\mathcal{O}(D,\mathcal{V})$ of sections.

The condition \eqref{eqn:1.2} 
 with regard to holomorphic actions on complex manifolds
has become the motivation that we introduced
 the concept of \textit{visible actions}  in \cite{xkglvis} .  
Recently, classification results on visible actions on complex manifolds
have been obtained in various settings,
see \cite{visiblesymm, xkgencar, Sa1,Sa2,Sa3}.
In this article,
we use a variant of visible actions,
namely, 
\textit{$S$-visible actions}; 
see Definition~\ref{def:Scompat}.

The second form of our 
multiplicity-free theorem is formalized as Theorem~\ref{thm:mfvis}
in terms of $S$-visible actions. 
Here, $S$ is a slice of the $H$-action on the base
space $D$.  
An old theorem of S. Kobayashi \cite{xkobasho}
 (see also Wolf \cite{xwolf})
 may be thought of as a {\it{propagation theorem 
 of irreducibility}} from fibers to the space of sections 
 when $S$ is a singleton.
 
The third form of our multiplicity-free theorem is formalized 
in the setting where the bundle $\mathcal{V} \to D$ is associated to a
principal bundle $K \to P \to D$ and to a representation
$(\mu, V)$ of the structure group $K$.
This is Theorem~\ref{thm:3.3}.
This formulation is useful for actual applications, in particular,
for branching problems
(decompositions of irreducible representations when restricted to
subgroups). 
In fact, this is the form that was used as a main machinery 
 of \cite{xkglvis,xkrims40} in finding various multiplicity-free theorems
 in concrete settings,
whereas the complete proof of Theorem \ref{thm:3.3}
(stated as \cite[Theorem 1.3]{xkglvis} and \cite[Theorem 2]{xkrims40}
loc.cit.)
has been postponed until the present article.

\medskip

\noindent
\textbf{Acknowledgement.}
A primitive case of Theorem \ref{thm:3.3} (the line bundle case)
together with its application to branching problems for semisimple
symmetric pairs $(G,H)$ was announced in \cite{xkmfjp}.
The heart of the proof of Theorem \ref{thm:2.2} is based on
reproducing kernels,
and was inspired by the original idea of Faraut--Thomas \cite{xft}.
I thank J. Faraut for enlightening discussions,
in particular,
for explaining the idea of \cite{xft} in an early stage of this work. 
Substantial part of its generalization in the present form was
obtained  during my visit to Harvard
University in 2000--2001.  
I express my gratitude to W. Schmid who gave me a wonderful
atmosphere of research.
M. Duflo suggested me to use the terminology ``propagation'' for
Theorem \ref{thm:main}.
Concrete applications of Theorem \ref{thm:main} and the theory of visible actions
 were presented in various occasions
including the Oberwolfach workshops organized by A. Huckleberry,
K.-H. Neeb and J. Wolf in 2000 and 2004
and at Winter School at Czech Republic organized by V. Sou\v cek in 2010.
A detailed account of the material of the present article
 (the proof of the propagation theorem) 
 was given in the
graduate course lectures at Harvard University (2008, spring
semester)
and at the University of Tokyo (2008, fall semester),
and also in a series of lectures at Functional Analysis X in Croatia (2008, summer).
I express my deep gratitude to the organizers and to the participants
for helpful and stimulating comments on various occasions.
Special thanks are due to Ms. Suenaga for her help in preparing for the final 
manuscript.

\numberwithin{equation}{subsection}

\section{%
Complex geometry and multiplicity-free theorem
}
This section gives a first form of our multiplicity-free theorem. 
We may regard it as a propagation theorem of multiplicity-free
property from fibers to spaces of sections in the setting where
there may exist infinitely many orbits on base spaces. 
The main result of this section is Theorem~\ref{thm:2.2}. 
We shall reformulate it by means of visible actions in
Section~\ref{sect:4},
and furthermore present its group theoretic version in Section~\ref{sect:3}.

\subsection{Equivariant holomorphic vector bundle}
\label{sect:2.1}
Let $\mathcal{V} = \amalg_{x\in D} \mathcal{V}_x \to D$ 
be a Hermitian holomorphic vector bundle
 over a connected complex manifold $D$.
We denote by $\OV$
 the space of holomorphic sections of $\mathcal{V} \to D$.
It carries a Fr\'echet topology
 by the uniform convergence on compact sets.

Suppose a Lie group $H$ acts on the bundle 
$\mathcal{V} \to D$
by automorphisms.
This means that
 the action of $H$ on the total space, 
denoted by $L_h : \mathcal{V} \to \mathcal{V}$,
 and the action
 on the base space, denoted simply by $h : D \to D, \ x \mapsto h \cdot x$,
 are both biholomorphic for $h \in H$,
 and that the induced linear map
 $L_h : \mathcal{V}_x \to \mathcal{V}_{h \cdot x}$ 
 between the fibers  is unitary for any $x \in D$.
In particular,
we have a unitary representation of the isotropy subgroup
$H_x := \{ h \in H : h \cdot x = x \}$
on the fiber $\mathcal{V}_x$.

The action of $H$ on the bundle $\mathcal{V} \to D$
 gives rise to a continuous representation on $\OV$
 by the pull-back of sections, namely, 
$s \mapsto L_h s(h^{-1} \cdot \ )$
 for $h \in H$ and $s \in \OV$.

\begin{defin}
\label{def:realizeO}
Suppose $\pi$ is a unitary 
 representation of $H$ defined on a Hilbert space
$\mathcal{H}$. 
We will say $\pi$ 
 is {\it realized in} $\OV$
 if there is an injective
 continuous $H$-intertwining map from $\mathcal{H}$ into $\OV$.
\end{defin}

Let $\{U_\alpha\}$ be trivializing neighborhoods of $D$,
 and $g_{\alpha \beta}: U_\alpha \cap U_\beta \to GL(n,\mathbb{C})$
 be the transition functions for the holomorphic vector bundle $\mathcal{V} \to D$.
Then, the anti-holomorphic vector bundle
 $\overline{\mathcal{V}} \to D$ is defined to be 
the complex vector bundle with
 the transition functions $\overline{g_{\alpha \beta}}$.
We denote by $\barOV$ 
 the space of anti-holomorphic sections for $\overline{\mathcal{V}} \to D$.

Suppose $\sigma$ is an anti-holomorphic diffeomorphism of $D$.
Then the pull-back $\sigma^* \mathcal{V} = \amalg_{x \in D} \mathcal{V}_{\sigma(x)}$
 is an anti-holomorphic vector
 bundle over $D$.
In turn,
 $\overline{\sigma^* \mathcal{V}} \to D$ is a holomorphic vector bundle 
over $D$.
The fiber at $x \in D$ is identified with
$\overline{\mathcal{V}_{\sigma(x)}}$,
the complex conjugate vector space of
$\mathcal{V}_{\sigma(x)}$
(see Section~\ref{sec:2.5.1}).

The holomorphic vector bundle $\overline{\sigma^* \mathcal{V}}$ is
isomorphic to $\mathcal{V}$ if and only if $\sigma$ lifts to an
anti-holomorphic endomorphism $\tilde{\sigma}$ of $\mathcal{V}$.
In fact, such $\tilde{\sigma}$ induces a conjugate linear isomorphism 
$\tilde{\sigma}_x: \mathcal{V}_x \to \mathcal{V}_{\sigma(x)}$,
which then defines a $\mathbb{C}$-linear isomorphism
\begin{equation}
\label{eqn:Psigma}
\Psi_x : \mathcal{V}_x \to (\overline{\sigma^* \mathcal{V}})_x,
\quad v \mapsto \overline{\tilde{\sigma}_x(v)}
\end{equation}
via the identification
$(\overline{\sigma^* \mathcal{V}})_x
 \simeq \overline{\mathcal{V}_{\sigma(x)}}$.
Then,
$\Psi : \mathcal{V} \to \overline{\sigma^* \mathcal{V}}$
is an isomorphism of holomorphic vector bundles such that its
restriction to the base space $D$ is the identity.
For simplicity, 
we shall use the letter $\sigma$ in place of $\tilde{\sigma}$.
For a Hermitian vector bundle $\mathcal{V}$,
by a bundle endomorphism $\sigma$,
we mean that $\sigma_x$ is furthermore 
isometric (or equivalently, $\Psi_x$ is unitary)
for any $x \in D$.

\subsection{Multiplicity-free theorem (first form)}
The following is a first form of our multiplicity-free theorem:
\begin{thm}\label{thm:2.2}  
Let $\mathcal{V} \to D$ be a Hermitian holomorphic vector bundle,
 on which a Lie group $H$ acts by automorphisms.
Assume:
\begin{texteqn}
\label{thmassB}
 the isotropy representation of $H_x$
 on the fiber $\mathcal{V}_x$ is multiplicity-free
 for any $x \in D$.
\end{texteqn}
We write its irreducible decomposition as
      $\mathcal{V}_x = \bigoplus\limits_{i=1}^{n(x)} \mathcal{V}_x^{(i)}$.
Assume furthermore that there exists an anti-holomorphic bundle endomorphism
      $\sigma$ 
satisfying the following two conditions:
for any $x \in D$, 
\begin{texteqn}
\label{thmassC}
there exists $h \in H$ such that 
$\sigma (x) = h \cdot x$, and
\end{texteqn}
\begin{texteqn}
\label{thmassA}
$\sigma_x(\mathcal{V}_x^{(i)}) = 
 L_h(\mathcal{V}_x^{(i)})$
for any $i$ $(1 \le i \le n(x))$.
\end{texteqn}

Then,
 any unitary representation that is 
 realized in $\OV$ is multiplicity-free.
\end{thm} 

We shall give a proof of Theorem~\ref{thm:2.2} in
Section~\ref{sect:2.5}. 

\begin{rem}\label{rem:2.2.2} 
\ \rm{ 1)}
The conditions \eqref{thmassB} -- \eqref{thmassA} 
of Theorem~\ref{thm:2.2} is local in the sense that
the same conclusion holds if $D'$ is an 
$H$-invariant open subset of $D$,
and if the conditions \eqref{thmassB} -- \eqref{thmassA} are satisfied
for $x \in D'$.
This is clear because the restriction map
$
\mathcal{O}(D,\mathcal{V}) \to \mathcal{O}(D',\mathcal{V} |_{D'})
$
is injective and continuous.  
\itm{2)}
The proof in Section~\ref{sect:2.5} shows that one
 can replace $H_x$ with its arbitrary subgroup $H_x'$
 in (\ref{thmassB}).  
(Such a replacement 
 makes (\ref{thmassB})  
stronger,
 and (\ref{thmassA}) 
weaker.)
\end{rem}  

In the following two subsections,
 we explain special cases of Theorem~\ref{thm:2.2}.

\subsection{Line bundle case}
We begin with the observation
 that the assumptions (\ref{thmassB})  
and (\ref{thmassA})  
are
 automatically fulfilled
 if $\mathcal{V}_x$ is irreducible, 
in particular,  
 if $\mathcal{V} \to D$ is a line bundle.
Hence, we have:
\begin{crl} 
In the setting of 
Theorem~\ref{thm:2.2},
assume 
 $\mathcal{V} \to D$ is a line bundle.
If 
there exists an anti-holomorphic bundle endomorphism 
satisfying \eqref{thmassC},
 then any unitary 
 representation that is 
realized in $\mathcal{O}(D,\mathcal{V})$ 
is multiplicity-free.
\end{crl}  
This case was announced in \cite{xkmfjp}\ with a sketch of proof, 
and its applications are extensively discussed in \cite{xkmfkorea}
for the branching problems
(i.e.\ the decomposition of the restriction of unitary
representations) 
 with respect to
reductive symmetric pairs.

\subsection{Trivial bundle case}
If the vector bundle is the trivial line bundle
 $\mathcal{V} = D \times \mathbb{C}$,
 then any anti-holomorphic diffeomorphism on $D$ lifts to
 an anti-holomorphic endomorphism of $\mathcal{V}$ by
$(x,u) \mapsto (\sigma(x), \bar{u})$.
Hence, we have:
\begin{crl}  
If there exists an anti-holomorphic diffeomorphism $\sigma$ of $D$
satisfying \eqref{thmassC}, 
then any unitary
 representation which is 
realized in $\mathcal{O}(D)$ is multiplicity-free.
\end{crl}  
This result was previously proved
  in Faraut and Thomas \cite{xft} under the assumption that 
$\sigma^2 = \operatorname{id}$.

\subsection{Propagation of irreducibility}
The strongest condition on group actions is transitivity.
Transitivity on base spaces
 guarantees that even irreducibility propagates from
fibers to spaces of sections.
The following result is due to S. Kobayashi \cite{xkobasho}. 

\begin{prop}
\label{prop:kobasho}
In the setting of Theorem~\ref{thm:2.2},
suppose that $H$ acts transitively on $D$ and that $H_x$ acts
irreducibly on $\mathcal{V}_x$ for some (equivalently, for any)
$x \in D$.
Then, there exists at most one unitary representation $\pi$ that can
be realized in $\mathcal{O}(D,\mathcal{V})$.
In particular, such 
$\pi$ is irreducible if exists.
\end{prop}

\begin{proof}
This is an immediate consequence of Lemma~\ref{lem:2.7} and
Proposition~\ref{prop:diagK} ($n(x)=1$ case) below,
which will be used in the proof of Theorem~\ref{thm:2.2} in
Section~\ref{sect:2.5}. 
\end{proof}

We note that the condition \eqref{thmassC}
is much weaker than the
 transitivity of the action of the group $H$ on $D$.
Our geometric condition \eqref{thmassC} brings us to 
the concept of visible actions,
which we shall discuss in Section~\ref{sect:4}.

\section{Proof of Theorem~\protect\ref{thm:2.2}}
\label{sect:2.5}
This section is devoted entirely to
 the proof of Theorem~\ref{thm:2.2}.   

\subsection{Some linear algebra}
\label{sec:2.5.1}
We begin carefully with basic notations. 

Given a complex Hermitian vector space $V$,
 we define a complex Hermitian vector space $\overline{V}$
 as a collection of the symbol $\overline{v}$ ($v \in V$)
 equipped with a scalar multiplication
 $a \bar{v} := \overline{\overline{a} v}$ for $a \in \mathbb{C}$,
 and with an inner product $(\bar{u}, \bar{v}) := (v,u)$.
\par
The complex dual space $\contra{V}$ is identified with $\overline{V}$
 by
 $\overline{V} \rarrowsim \contra{V}, \ \bar{v} \mapsto (\cdot , v)$.
In particular,
 we have a natural isomorphism of complex vector spaces:
\begin{equation}\label{eqn:2.5.1}
   V \otimes \overline{V} \rarrowsim \End(V).
\end{equation}
\par
Given a unitary map $A : V \to W$ between Hermitian vector spaces,
 we define a unitary map
$
   \overline{A}:\overline{V} \to \overline{W}
$
by
$
                \overline{v}\mapsto  \overline{Av}
$.
Then the induced map 
$
   A \otimes \overline {A}:V \otimes \overline{V} \to W \otimes \overline{W}
$
 gives rise to a complex linear isomorphism: 
\begin{equation}\label{eqn:2.5.2}
    A_\sharp : \End(V) \to \End(W).
\end{equation}
Then, it is readily seen from the unitarity of $A$ that
\begin{equation}\label{eqn:2.5.4}
    A_\sharp   (\operatorname{id}_V)
  =
                \operatorname{id}_W.
\end{equation}
In particular, if an endomorphism of $V$ is diagonalized with respect to
 an orthogonal direct sum decomposition
 $V = \bigoplus_{i=1}^n V^{(i)}$, then we have the following formula
 of 
$A_{\sharp}$: 
\begin{equation}\label{eqn:2.5.5}
     A_\sharp \left(\sum_{i=1}^n \lambda_i \operatorname{id}_{V^{(i)}} \right)
  =
     \sum_{i=1}^n \lambda_i \operatorname{id}_{A ({V^{(i)}})}
 \quad (\lambda_1,\dots,\lambda_n \in \mathbb{C}).
\end{equation}

\subsection{Reproducing kernel for vector bundles}
This subsection summarizes some basic results on reproducing kernels
for holomorphic vector bundles.
The results here are  standard  for the trivial
bundle case.

Suppose we are given 
 a continuous embedding 
$\mathcal{H} \hookrightarrow \mathcal{O}(D,\mathcal{V})$
of a Hilbert space $\mathcal{H}$ into the
Fr\'{e}chet space $\mathcal{O}(D,\mathcal{V})$ of holomorphic sections
of the holomorphic vector bundle $\mathcal{V} \to D$.
The continuity implies in particular that for each $y \in D$ 
the point evaluation map: 
$$
\mathcal{H} \to \mathcal{V}_y,
\quad
f \mapsto f(y)
$$
is continuous.
Then, by the Riesz representation theorem,
there exists uniquely
$K_{\mathcal{H}} (\cdot, y) \in \mathcal{H} \otimes
\overline{\mathcal{V}_y} $
such that
\begin{equation}\label{eqn:repro}
(f, K_{\mathcal{H}} (\cdot, y) )_{\mathcal{H}} = f(y)
\quad \text{for any $f \in \mathcal{H}$.}
\end{equation}

We take an orthonormal basis $\{ \varphi_\nu \}$ of $\mathcal{H}$,
and expand $K_{\mathcal{H}}$ as
\begin{equation}
\label{eqn:Kexpan}
K_{\mathcal{H}} (\cdot, y) = 
    \sum_{\nu} a_{\nu} (y) \varphi_{\nu} (\cdot). 
\end{equation}
It follows from \eqref{eqn:repro} that
 the coefficient $a_{\nu} (y)$ is given by
$$
a_{\nu} (y) = (K_{\mathcal{H}} (\cdot, y), 
\varphi_{\nu} (\cdot) )_{\mathcal{H}}
= \overline{\varphi_{\nu} (y)},
$$
and the expansion of $K_\mathcal{H}$ 
converges in $\mathcal{H}$. 
By the continuity 
$\mathcal{H}\hookrightarrow\mathcal{O}(D,\mathcal{V})$ 
again,
\eqref{eqn:Kexpan} converges uniformly on each compact set 
for each fixed $y \in D$.
Thus,
 $K_{\mathcal{H}}(x,y)$
 is given by the formula: 
\begin{equation}\label{eqn:2.6.1}
   K_{\mathcal{H}}(x,y) \equiv  K(x,y)
  = \sum_{\nu} \varphi_\nu(x) \overline{\varphi_\nu(y)}
          \in \mathcal{V}_x \otimes \overline{\mathcal{V}_y} \; ,
\end{equation}
and defines 
a smooth section of the exterior tensor product bundle 
$\mathcal{V} \boxtimes \overline{\mathcal{V}}
 \to D \times D$
which is holomorphic in the first argument and anti-holomorphic in the
second. 
We will say $K_{\mathcal{H}}$ is the \textit{reproducing kernel} of the
Hilbert space $\mathcal{H} \subset \mathcal{O}(D,\mathcal{V})$. 

For the convenience of the reader,
we pin down basic properties of reproducing kernels for holomorphic
vector bundles in a way that we use later.

\begin{lem}\label{lem:2.6}  
{\rm 1)}\enspace
Let $K_i(x,y)$ 
 be the reproducing kernels of
 Hilbert spaces $\mathcal{H}_i \subset \OV$
 with inner products $(\ , \ )_{\mathcal{H}_i}$,
 respectively, for $i = 1, 2$.
If $K_1 \equiv K_2$,
 then the two subspaces $\mathcal{H}_1$ and $\mathcal{H}_2$ coincide
and the inner products $(\ , \ )_{\mathcal{H}_1}$ and
$(\ , \ )_{\mathcal{H}_2}$ are the same.
\itm{2)}
If $K_1(x,x) = K_2(x,x)$ for all $x \in D$,
 then $K_1 \equiv K_2$.
\end{lem}  
\begin{proof}
{\rm 1)}\enspace
Let us reconstruct the Hilbert space $\mathcal{H}$
from the reproducing kernel $K$.
For each $y \in D$ and $v^* \in \mathcal{V}^*_y :=
\overline{\mathcal{V}_y}^\vee
$,
 we define $\psi(y,v^*)$ by
$$
 \psi(y, v^*) := \langle K(\cdot, y),v^* \rangle  \in 
 \mathcal{H}.
$$
Here, $\langle \ , \ \rangle$ denotes the canonical pairing between
$\overline{\mathcal{V}_y}$ and $\overline{\mathcal{V}_y}^\vee$.
We claim that 
 the $\mathbb{C}$-span of 
$\set{\psi(y,v^*)}{y \in D, v^* \in \mathcal{V}^*_y}$ 
is dense in $\mathcal{H}$. 
This is because 
$(f, \psi(y,v^*))_{\mathcal{H}} = \langle f(y),v^* \rangle$
for any $f \in \mathcal{H}$ by \eqref{eqn:repro}.
Thus, the Hilbert space $\mathcal{H}$ is reconstructed from the
  pre-Hilbert structure
\begin{equation}
  (\psi(y_1, v^*_1), \psi(y_2, v^*_2))_{\mathcal{H}}
  := \langle K(y_2, y_1), v^*_2 \otimes \overline{v^*_1} \rangle.
\end{equation}
\itm{2)}
We denote by $\overline{D}$ the complex manifold
 endowed with the conjugate complex structure on the same real
manifold $D$.
Then,
 $\overline{\mathcal{V}} \to \overline{D}$ is a holomorphic vector bundle,
 and we can form a holomorphic vector bundle
 $\mathcal{V} \boxtimes \overline{\mathcal{V}} \to D \times \overline{D}$.
In turn, $K_i(\cdot,\cdot)$ are regarded as
its holomorphic sections. 
As the diagonal embedding
 $\iota : D \to D \times \overline{D}, z \mapsto (z, z)$
 is totally real,
our assumption 
 $(K_1-K_2)|_{\iota(D)} \equiv 0$ implies 
 $K_1-K_2\equiv 0$ by the unicity theorem of holomorphic functions.
\end{proof}  

\subsection{Equivariance of the reproducing kernel}
Next, suppose that the Hermitian holomorphic vector bundle 
$\mathcal{V} \to D$ is $H$-equivariant and that
$(\pi, \mathcal{H})$ is a unitary representation of $H$
 realized in $\OV$.
Let $K_{\mathcal{H}}$ be the reproducing kernel of the embedding
$\mathcal{H} \hookrightarrow \mathcal{O}(D,\mathcal{V})$.
We shall see how the unitarity of $(\pi,\mathcal{H})$ is reflected in
the reproducing kernel $K_{\mathcal{H}}$.

We regard 
$K_\mathcal{H}(x,x) \in \mathcal{V}_x \otimes \overline{\mathcal{V}_x}$
  as an element of $\End(\mathcal{V}_x)$ via the isomorphism 
\eqref{eqn:2.5.1}.   
Then, we have: 
\begin{lem}\label{lem:2.7}  
With the notation \eqref{eqn:2.5.2}  
applied to $L_h: \mathcal{V}_x \to \mathcal{V}_{h \cdot x}$,
 we have 
$$
    K_\mathcal{H}(h \cdot x, h \cdot x) =  (L_h)_\sharp K_\mathcal{H}(x,  x)
    \quad
     \text{for any $h \in H$}.  
$$
In particular, 
  $K_\mathcal{H}(x,x) \in \End_{H_x}(\mathcal{V}_x)$
for any $x \in D$.   
 \end{lem}  
\begin{proof}  
Let $\{\varphi_\nu\}$ be an orthonormal basis of $\mathcal{H}$.
Since $(\pi, \mathcal{H})$ is a unitary representation,
 $\{\pi(h)^{-1} \varphi_\nu\}$ is also an orthonormal basis of $\mathcal{H}$
 for every fixed $h \in H$.
Because the formula \eqref{eqn:2.6.1} of
the reproducing kernel is valid for any orthonormal basis, 
 we have 
\begin{align}  
   K_{\mathcal{H}}(x, y)
 &= \sum_\nu (\pi(h)^{-1} \varphi_\nu)(x) 
             \overline{(\pi(h)^{-1} \varphi_\nu)(y)}
\\
 & = \sum_\nu L_{h^{-1}} \varphi_\nu(h \cdot x) 
              \overline{L_{h^{-1}} \varphi_\nu(h \cdot y)}
\nonumber
\\
 & = (L_{h^{-1}} \otimes \overline{L_{h^{-1}}}) 
      K_\mathcal{H}(h \cdot x,  h \cdot y) 
\nonumber
\end{align}  
for any $x, y \in D$.
Hence,  
 $(L_{h} \otimes \overline{L_{h}}) K_{\mathcal{H}}(x, y) = K_\mathcal{H}(h \cdot x,  h \cdot y)$
and  Lemma follows.
\end{proof}  

\subsection{Diagonalization of the reproducing kernel}
The reproducing kernel for a holomorphic vector bundle is a matrix
valued section as we have defined in \eqref{eqn:2.6.1}.
The multiplicity-free property of the isotropy representation 
on the fiber diagonalizes 
the reproducing kernel:

\begin{prop}\label{prop:diagK}
Suppose $(\pi,\mathcal{H})$ is a unitary representation of $H$
realized in $\mathcal{O}(D,\mathcal{V})$.
Assume that the isotropy representation of $H_x$ on the fiber
$\mathcal{V}_x$ decomposes as a multiplicity-free sum of 
 irreducible representations of $H_x$ as 
$\mathcal{V}_x = \bigoplus\limits_{i=1}^n \mathcal{V}_x^{(i)}$. 
(Here, $n \equiv n(x)$ may depend on $x \in D$.) 
Then, the reproducing kernel is of the form 
$$
K_{\mathcal{H}} (x,x) = \sum_{i=1}^n \lambda^{(i)} (x)
   \operatorname{id}_{\mathcal{V}_{x}^{(i)}}
$$
for some complex numbers
$\lambda^{(1)}(x),\dots,\lambda^{(n)}(x)$.
\end{prop}

\begin{proof}
A direct consequence of Lemma~\ref{lem:2.7} and Schur's lemma.
\end{proof}

\subsection{Construction of an anti-linear isometry $J$}
In the setting of 
Theorem~\ref{thm:2.2}, 
suppose that $\sigma$ is an anti-holomorphic bundle endomorphism. 
We define a conjugate linear map
\begin{equation}
J : \OV \to \OV,
 \quad f\mapsto \sigma^{-1} \circ f \circ \sigma,
\end{equation}
namely,
 $J f(x) := \sigma^{-1}(f(\sigma(x)) $
for $x \in D$.  

\begin{lem}\label{lem:2.8}  
If the conditions \eqref{thmassB} -- \eqref{thmassA} are satisfied,
then 
 $J$ is an isometry from $\mathcal{H}$ onto $\mathcal{H}$
for any unitary representation $(\pi,\mathcal{H})$
realized in $\mathcal{O}(D,\mathcal{V})$. 
\end{lem}  
\begin{proof}  
We define a Hilbert space $\widetilde{\mathcal{H}} := J(\mathcal{H})$,
 equipped with the inner product
\begin{equation*}
  (J f_1, J f_2)_{\widetilde{\mathcal{H}}} := (f_2, f_1)_{\mathcal{H}}
 \quad \text{for } f_1, f_2 \in \mathcal{H}.
\end{equation*}
Let us show that
 the reproducing kernel $K_{\widetilde{\mathcal{H}}}$
for $\widetilde{\mathcal{H}}$
coincides with $K_{\mathcal{H}}$.
To see this, we take
 an orthonormal basis  $\{\varphi_\nu\}$ of $\mathcal{H}$.
 Then, 
 $\{J \varphi_\nu\}$ is an orthonormal basis of
 $\widetilde{\mathcal{H}}$, and therefore
\begin{align*}  
 K_{\widetilde{\mathcal{H}}}(x,y) 
 &=\sum_\nu J \varphi_\nu(x) \overline{J \varphi_\nu(y)}
\\
 &=\sum_\nu \sigma_x^{-1}\left( {\varphi_\nu(\sigma(x))}\right) \ 
 \overline{\sigma_y^{-1}\left({\varphi_\nu(\sigma(y))}\right)}
\\
 &= \left({\sigma_x^{-1}} \otimes \overline{\sigma_y^{-1}}\right)
        K_{\mathcal{H}}(\sigma(x), \sigma(y)).
\end{align*}  
For $x = y$, 
this formula can be restated as
\begin{equation}
   K_{\widetilde{\mathcal{H}}}(x,x)  
   =
   (\sigma_x^{-1})_\sharp
        K_{\mathcal{H}}(\sigma(x), \sigma(x))
\end{equation}
 with the notation \eqref{eqn:2.5.2}  
applied to the unitary map 
$\sigma_x^{-1}:{\mathcal{V}_{\sigma(x)}}\to \mathcal{V}_x$. 
We fix $x \in D$,
and take $h \in H$ such that $\sigma (x) = h \cdot x$
 (see \eqref{thmassC}). 
Then,
\begin{equation}\label{eqn:2.8.3}
 K_{\widetilde{\mathcal{H}}}(x,x) 
= (\sigma_x^{-1})_\sharp
        K_{\mathcal{H}}(h \cdot x, h \cdot x)
=  (\sigma_x^{-1})_\sharp (L_h)_\sharp K_{\mathcal{H}}(x, x).
\end{equation}
Here,
 the last equality follows from Lemma~\ref{lem:2.7}.  
 
Since the action of $H_x$ on $\mathcal{V}_x$ is multiplicity-free, 
 it follows from Proposition~\ref{prop:diagK}   
that
 there exist complex numbers $\lambda^{(i)}(x)$ such that
\begin{equation*}
     K_\mathcal{H}(x,x) = \sum_i \lambda^{(i)}(x) \operatorname{id}_{\mathcal{V}_x^{(i)}}.
\end{equation*}
Then by 
\eqref{eqn:2.5.4}  
we have
\begin{equation}\label{eqn:2.8.4}
  (L_h)_\sharp K_\mathcal{H}(x,x) = \sum_i \lambda^{(i)}(x) \operatorname{id}_{L_h(\mathcal{V}_x^{(i)})}.
\end{equation}
Furthermore,
 since $\sigma_x^{-1}({L_h(\mathcal{V}_x^{(i)})}) = \mathcal{V}_x^{(i)}$ 
 by the assumption (\ref{thmassA}), 
 it follows from 
\eqref{eqn:2.5.5}  
that
\begin{equation}\label{eqn:2.8.5}
 (\sigma_x^{-1})_\sharp \left(\sum_i \lambda^{(i)}(x) \operatorname{id}_{L_h(\mathcal{V}_x^{(i)})}\right)
 =
  \sum_i \lambda^{(i)} (x)\operatorname{id}_{\mathcal{V}_x^{(i)}}.
\end{equation}
Combining \eqref{eqn:2.8.3},  
\eqref{eqn:2.8.4}  
and 
\eqref{eqn:2.8.5},  
we get
$$
    K_{\widetilde{\mathcal{H}}}(x,x) = K_\mathcal{H}(x,x).
$$
Then,
 by Lemma~\ref{lem:2.6},  
 the Hilbert space $\widetilde{\mathcal{H}}$ coincides with $\mathcal{H}$
  and 
$$
  (J f_1, J f_2)_{\mathcal{H}} = (Jf_1, Jf_2)_{\widetilde{\mathcal{H}}}
   = (f_2, f_1)_{\mathcal{H}}
 \quad \text{for } f_1, f_2 \in \mathcal{H}.  
$$
This is what we wanted to prove.
\end{proof}  

\begin{rem}
In terms of the bundle isomorphism
$\Psi: \mathcal{V} \to \overline{\sigma^*\mathcal{V}}$
(see \eqref{eqn:Psigma}),
$J$ is given by
$(Jf)(x) = \Psi_x^{-1}(\overline{f(\sigma(x))})$. 
We note
$$
   J^2 = \operatorname{id}
    \quad
    \text{on } \OV
$$
if $\sigma^2 = \operatorname{id}_{\mathcal{V}}$,
or equivalently, if
$\sigma^2 = \operatorname{id}_D$ and 
$\overline{\Psi_{\sigma(x)}}\circ \Psi_x = \operatorname{id}_{\mathcal{V}_x}$
 for any $x \in D$.
However, we do not use this condition to prove Theorem~\ref{thm:2.2}.
\end{rem}

\subsection{Proof of Theorem~\protect\ref{thm:2.2}}
As a final step,
we need the following lemma which was proved in \cite{xft}
under the assumption that $J^2 = \operatorname{id}$
and that $\mathcal{V} \to D$ is the trivial line bundle.
For the sake of completeness,
 we give a proof here.
\begin{lem}\label{lem:2.9}  
For $A \in \End_H(\mathcal{H})$,
 the adjoint operator $A^*$ is given by
\begin{equation}
    A^* = J A J^{-1}.
\end{equation}
\end{lem}  
\begin{proof}  
We divide the proof into two steps.
\newline
\textbf{Step 1} (self-adjoint case).\enspace
We may and do assume that $A-I$ is positive definite
 because neither the assumption nor the conclusion changes
 if we replace $A$ by $A + c I$ $(c \in \mathbb{R})$.
 Here, we note that $A + c I$ is positive definite if $c$
is greater than the operator norm $\| A \|$.

{}From now,
assume $A \in \End_H(\mathcal{H})$ is a 
 self-adjoint operator such that $A-I$ is positive definite.
We introduce a 
 pre-Hilbert structure on $\mathcal{H}$ by
\begin{equation}
  (f_1, f_2)_{\mathcal{H}_A} := (A f_1, f_2)_{\mathcal{H}}
 \quad \text{for } f_1, f_2 \in \mathcal{H}.  
\end{equation}
Since $A-I$ is positive definite,
we have
$$
(f,f)_{\mathcal{H}} \le (f,f)_{\mathcal{H}_A}
   \le \| A \| (f,f)_{\mathcal{H}}
\quad \text{for $f \in \mathcal{H}$.}
$$
Therefore,
$\mathcal{H}$ is still complete with respect to the new inner product
$(\ , \ )_{\mathcal{H}_A}$.
The resulting Hilbert space will be denoted by $\mathcal{H}_A$.

If $f_1, f_2 \in \mathcal{H}$ and $g \in H$,
 then
\begin{align*}  
   (\pi(g) f_1, \pi(g) f_2)_{{\mathcal{H}}_A}
  &= (A \pi(g) f_1, \pi(g) f_2)_{\mathcal{H}}
\\
  &= (\pi(g) A f_1, \pi(g) f_2)_{\mathcal{H}}
  = (A f_1,  f_2)_{\mathcal{H}}
  = (f_1,  f_2)_{{\mathcal{H}}_A}.
\end{align*}  
Therefore,
 $\pi$ also defines a unitary representation on ${\mathcal{H}}_A$.
Applying Lemma~\ref{lem:2.8}   
to both $\mathcal{H}_A$ and $\mathcal{H}$,
 we have 
\begin{multline*}  
    (A f_1, f_2)_{\mathcal{H}}
 =  (f_1, f_2)_{\mathcal{H}_A}
 =  (J f_2, J f_1)_{\mathcal{H}_A}
 =  (A J f_2, J f_1)_{\mathcal{H}}
\\
 =  (J f_2, A^* J f_1)_{\mathcal{H}}
 =  (J f_2, J J^{-1} A^* J f_1)_{\mathcal{H}}
 =  (J^{-1} A^* J f_1, f_2)_{\mathcal{H}}.
\end{multline*}  
Hence, $A = J^{-1} A^* J$.

\noindent
\textbf{Step 2} (general case).\enspace
Suppose $A \in \End_H(\mathcal{H})$. 
Then $A^*$ also commutes with $\pi(g)$ ($g \in H$) because $\pi$ is unitary.
We put $B := \frac{1}{2}(A+ A^*)$ and $C := \frac{\sqrt{-1}}{2}(A^* - A)$.
Then,
 both $B$ and $C$ are self-adjoint operators commuting with $\pi(g)$
 ($g \in H$).
It follows from Step 1 that
$
    B^* = J B J^{-1}
$
and
$
    C^* = J C J^{-1}.
$
Since $J$ is conjugate-linear,
 we have
$   (\sqrt{-1} C)^*  = J (\sqrt{-1} C) J^{-1}$.
Hence,
 $A = B + \sqrt{-1} C$ also satisfies
$
    A^* = J A J^{-1}.
$
\end{proof}  

\begin{proof}[Proof of Theorem~\ref{thm:2.2}.] 
Let $A$, $B \in \End_H(\mathcal{H})$. 
By Lemma~\ref{lem:2.9},  
 we have
$$
   A B = J^{-1} (A B)^* J = J^{-1} B^* J J^{-1} A^* J = B A. 
$$
Therefore,
the ring $\End_H(\mathcal{H})$  is commutative.
\end{proof}

\section{Visible actions on complex manifolds}
\label{sect:4}

This section analyzes the geometric condition \eqref{thmassC} on the
complex manifold $D$.
We shall introduce the concept of $S$-visible actions,
with which Theorem~\ref{thm:2.2} is reformulated in a simpler manner
(see Theorem~\ref{thm:mfvis}).

\subsection{Visible actions on complex manifolds}
\label{subsec:visible}
Suppose a Lie group $H$ acts holomorphically on a connected complex
manifold $D$.

\begin{defin}
\label{def:S}
We say the action is \textit{$S$-visible}  
if there exist 
a subset $S$ of $D$ such that
\begin{texteqn}
\label{eqn:S0} 
\text{$D' := H \cdot S $ is open in $D$,}
\end{texteqn}
and an anti-holomorphic diffeomorphism $\sigma$ of $D'$ 
satisfying the following two conditions:
\begin{texteqn}
\label{eqn:S2}
$\sigma|_S = \operatorname{id}$,
\end{texteqn}
\begin{texteqn}
\label{eqn:S3}
$\sigma$ preserves every $H$-orbit in $D'$.
\end{texteqn}
\end{defin}

\begin{rem}
\label{rem:Slocal}
The above condition is \textit{local} in the sense that we may replace
$S$ by its subset $S'$ in Definition~\ref{def:S} as far as
$H \cdot S'$ is open in $D$.
\end{rem}

\begin{rem}
\label{rem:S} 
By the definition of $D'$,
it is obvious that 
\begin{texteqn}
\label{eqn:S1}
$S$ meets every $H$-orbit in $D'$.
\end{texteqn}
Thus,
Definition~\ref{def:S} is essentially the same with 
\textit{strong visibility} in the sense of 
 \cite[Definition~3.3.1]{xkrims40}. 
In fact, the difference is only an additional requirement that $S$ is
a smooth submanifold in \cite{xkrims40}.
We note that
if $S$ is a smooth submanifold in Definition~\ref{def:S}, 
then $S$ is totally real by the condition
\eqref{eqn:S2}, and consequently, the $H$-action becomes
\textit{visible} in the sense of \cite{xkglvis}
(see \cite[Theorem~4.3]{xkrims40}). 
\end{rem}

\subsection{Compatible automorphism}
\label{subsec:Scompatible}
Retain the setting of Definition~\ref{def:S}. 
Suppose
$\sigma$ is the anti-holomorphic diffeomorphism of $D'$.
Twisting the original $H$-action by $\sigma$,
we can define another holomorphic action of $H$ on $D'$ by
$$
D' \to D', \ 
x \mapsto \sigma(h \cdot \sigma^{-1}(x)).
$$
If this action can be realized by $H$,
namely, if there exists
 a group automorphism $\tilde{\sigma}$ of $H$ such that 
$$
\tilde{\sigma}(h)\cdot x = \sigma(h \cdot \sigma^{-1}(x))
\quad\text{for any $x \in D'$,}
$$
we say $\tilde{\sigma}$ is \textit{compatible} with $\sigma$.
This condition is restated simply as 
\begin{equation}
\label{eqn:compatHD}
\tilde{\sigma}(h) \cdot \sigma(y) = \sigma(h \cdot y)
\quad\text{for any $y \in D'$.}
\end{equation}

\begin{defin}
\label{def:Scompat}
We say an $S$-visible action has a \textit{compatible
automorphism} of the transformation group $H$ if there exists 
 an automorphism $\tilde{\sigma}$ of the group $H$ satisfying
the condition \eqref{eqn:compatHD}. 
\end{defin}

We remark that the condition \eqref{eqn:S3} follows from
\eqref{eqn:S0} and \eqref{eqn:S2} if there exists $\tilde{\sigma}$
satisfying \eqref{eqn:compatHD}.
In fact, any $H$-orbit in $D'$ is of the form $H \cdot x$ for some 
$x \in S$, and then
$$
\sigma(H \cdot x) = \tilde{\sigma}(H) \cdot \sigma(x) = H \cdot x
$$
by \eqref{eqn:S2} and \eqref{eqn:compatHD}.

Suppose $\mathcal{V} \to D$ is an $H$-equivariant holomorphic vector
bundle. 
If there is a compatible automorphism $\tilde{\sigma}$
of $H$
with an anti-holomorphic diffeomorphism $\sigma$ on $D$, 
then we have the following isomorphism:
$$
( \overline{\sigma^* \mathcal{V}} ) _{h\cdot y}
   \simeq   \overline{\mathcal{V}}_{\sigma(h\cdot y)}
   =   \overline{\mathcal{V}}_{\tilde{\sigma}(h)\cdot\sigma(y)}
 \quad\text{for $h \in H$ and $y \in D$}.
$$
Therefore,
we can let $H$ act equivariantly on the holomorphic vector bundle 
$\overline{\sigma^* \mathcal{V}} \to D$
 by defining the left translation on 
 $\overline{\sigma^*\mathcal{V}}$ as
$$
L_h^\sigma : ( \overline{\sigma^*\mathcal{V}} ) _y
   \to
             ( \overline{\sigma^*\mathcal{V}} ) _{h\cdot y}
$$
via the identification with the left translation
$
\overline{L_{\tilde{\sigma}(h)}} :
    \overline{\mathcal{V}}_{\sigma(y)}
  \to
    \overline{\mathcal{V}}_{\tilde{\sigma}(h)\cdot\sigma(y)} .
$
Then, the two $H$-equivariant holomorphic vector bundles $\mathcal{V}$
and $\overline{\sigma^* \mathcal{V}}$ are isomorphic if and only if
$\sigma$ lifts to an anti-holomorphic bundle endomorphism 
$\sigma$ (we use the same letter) 
which respects the $H$-action 
in the sense that
\begin{equation}
\label{eqn:HisomV}
L_{\tilde{\sigma}(h)} \circ \sigma = \sigma \circ L_h
\quad\text{on $\mathcal{V}$ \quad for any $h \in H$.}
\end{equation}

\subsection{Propagation of multiplicity-free property}
By using the concept of $S$-visible actions,
we give a second form of our main theorem as follows:

\begin{thm}
\label{thm:mfvis}
Let $\mathcal{V}\to D$ be an $H$-equivariant Hermitian holomorphic
vector bundle.
Assume the following three conditions are satisfied:
\begin{texteqn}
\label{eqn:42B}
{\upshape(}Base space{\upshape)} \
The action on the base space $D$ is $S$-visible with a compatible
automorphism of the group $H$ (Definition~\ref{def:Scompat}).
\end{texteqn}
\begin{texteqn}
\label{eqn:42F}
{\upshape(}Fiber{\upshape)} \
The isotropy representation of $H_x$ on $\mathcal{V}_x$ is
multiplicity-free for any $x \in S$.
\end{texteqn}
We write its irreducible decomposition as
$$
\mathcal{V}_x = \bigoplus_{i=1}^{n(x)} \mathcal{V}_x^{(i)}.
$$
\begin{texteqn}
\label{eqn:42C}
{\upshape(}Compatibility{\upshape)} \
$\sigma$ lifts to an anti-holomorphic endomorphism
(we use the same letter $\sigma$) of the 
 $H$-equivariant Hermitian holomorphic vector bundle $\mathcal{V}$ such
that \eqref{eqn:HisomV} holds and
\end{texteqn}
\begin{equation}
\sigma_x (\mathcal{V}_x^{(i)}) = \mathcal{V}_x^{(i)}
\quad\text{for $1 \le i \le n(x)$, $x\in S$.}
\tag{\ref{eqn:42C})(a}
\end{equation}
Then, any unitary representation which is realized in
$\mathcal{O}(D,\mathcal{V})$
is multiplicity-free.
\end{thm}

The difference between the conditions of Theorem \ref{thm:mfvis} with
 the previous conditions 
\eqref{thmassB} and \eqref{thmassA}
in Theorem~\ref{thm:2.2} is the following:
The conditions
\eqref{eqn:42F} and
\eqref{eqn:42C}(a) are imposed only on the slice $S$, 
while the conditions in Theorem~\ref{thm:2.2} were imposed on the
 whole base space $D$
(or at least its open subset).

\begin{rem}
\label{rem:mfvis}
We can sometimes find a slice $S$ such that the isotropy subgroup
$H_x$ is independent of \textit{generic} $x \in S$.
Bearing this in mind,
we set 
\begin{align*}
H_S := {}
& \bigcap_{x \in S} H_x
\\
= {}
& \{ g \in H : gx = x \quad\text{for any $x \in S$} \}.
\end{align*}
Theorem~\ref{thm:mfvis} still holds if we replace $H_x$ 
with $H_S$ (see also Remark~\ref{rem:2.2.2}(2)).
\end{rem}

\begin{proof}
We shall reduce Theorem~\ref{thm:mfvis} to Theorem~\ref{thm:2.2} by
using  the $H$-equivariance of
the bundle endomorphism $\sigma$. 
Let us show that the conditions (\ref{thmassB}), (\ref{thmassC}) and
(\ref{thmassA}) are satisfied for the $H$-invariant open subset
$D' := H \cdot S$ of $D$. 

First we observe that the
 condition \eqref{eqn:S3} implies (\ref{thmassC}) because
$\sigma(x) \in \sigma(H \cdot x) = H \cdot x$
for any $x \in D'$.

Next, take any element $x \in D'$ and we write $x = h \cdot x_0$
($h \in H$, $x_0 \in S$).
We set
\begin{equation*}
\mathcal{V}_x^{(i)} 
:=
  L_h(\mathcal{V}_{x_0}^{(i)})
\quad (1 \le i \le n(x_0)). 
\end{equation*}
Through the group isomorphism
$H_{x_0} \rarrowsim H_{x}, l \mapsto hlh^{-1}$
and the left translation
$L_h: \mathcal{V}_{x_0} \to \mathcal{V}_{x}$, 
we get the isomorphism between the two isotropy representations,
$H_{x_0} \to GL(\mathcal{V}_{x_0})$
and
$H_{x} \to GL(\mathcal{V}_{x})$,
because
$L_{hlh^{-1}} = L_h \circ L_l \circ L_h^{-1}$
$(l \in H_{x_0})$.
In particular,
the direct sum
$$
\mathcal{V}_x=\bigoplus_{i=1}^{n(x_0)}\mathcal{V}_x^{(i)}
$$
gives a multiplicity-free decomposition of irreducible representations
of $H_x$.
Hence, the condition (\ref{thmassB}) is satisfied for all
$x \in D'$. 

Finally, we set $g := \tilde{\sigma}(h) h^{-1} \in H$.
As $\sigma(x_0) = x_0$, 
we have
$$
\sigma(x) = \sigma(h \cdot x_0) = \tilde{\sigma}(h) \cdot \sigma(x_0)
= \tilde{\sigma}(h) \cdot x_0 = g \cdot x.
$$
Besides, we have for any $i$ $(1 \le i \le n(x) = n(x_0))$,
\begin{align*}
\sigma_x(\mathcal{V}_x^{(i)})
&= \sigma_x (L_h (\mathcal{V}_{x_0}^{(i)}))
\\
&= L_{\tilde{\sigma}(h)} (\sigma_{x_0} (\mathcal{V}_{x_0}^{(i)}))
&&\text{by \eqref{eqn:HisomV}}
\\
&= L_{\tilde{\sigma}(h)} \left(\mathcal{V}_{x_0}^{(i)}\right)
&&\text{by \eqref{eqn:42C}(a)}
\\
&= L_{\tilde{\sigma}(h)h^{-1}} L_h (\mathcal{V}_{x_0}^{(i)})
\\
&= L_g(\mathcal{V}_x^{(i)}).
\end{align*}
Hence, the condition (\ref{thmassA}) holds for any
$x \in D'$. 
Therefore, all the assumptions of Theorem~\ref{thm:2.2} are satisfied
for the open subset $D'$.
Now,
Theorem~\ref{thm:mfvis} follows from 
Theorem~\ref{thm:2.2} and Remark~\ref{rem:2.2.2} (1). 
\end{proof}

\section{Multiplicity-free theorem for associated bundles}
\label{sect:3}
This section provides a third form of our multiplicity-free theorem 
(see Theorem~\ref{thm:3.3}). 
It is intended for actual applications to group representation theory,
especially to branching problems. 
The idea here is to reformalize the geometric condition of
Theorem~\ref{thm:mfvis} (second form) in terms of the representation
of the structure group of an equivariant principal bundle. 

Theorem~\ref{thm:3.3} is used as a main machinery in
 \cite{xkglvis,xkrims40}
(referred to as
 \cite[Theorem~1.3]{xkglvis} and \cite[Theorem~2]{xkrims40},
of which we have postponed the proof to this article)
for various multiplicity-free theorems including the
 following cases:
\begin{itemize}
\item
tensor product representations of $GL(n)$ \cite[Theorem~3.6]{xkglvis},
\item
branching problems for 
$GL(n)\downarrow GL(n_1)\times GL(n_2)\times GL(n_3)$ 
(\cite[Theorem~3.4]{xkglvis}),
\item
Plancherel formulae for vector bundles over Riemannian symmetric
spaces
(\cite[Theorems~21 and 30]{xkrims40}).
\end{itemize}

\subsection{Automorphisms on equivariant principal bundles}
\label{sect:3.1}

We begin with the setting where a Hermitian holomorphic vector
bundle $\mathcal{V}$ over a connected complex manifold $D$ 
is given as the associated bundle
$\mathcal{V} \simeq P \times_K V$
to the following data $(P,K,\mu,V)$:
\begin{align*}
&\text{$K$ is a Lie group,}
\\
&\text{$\varpi : P \to D$ is a principal $K$-bundle,}
\\
&\text{$V$ is a finite dimensional Hermitian vector space,}
\\
&\text{$\mu : K \to GL_{\mathbb{C}}(V)$ is a unitary representation.}
\end{align*}

Suppose that a Lie group $H$  acts on $P$ from the left,
 commuting with the right action of $K$.  
Then $H$ acts also on the Hermitian vector bundle
 $\mathcal{V} \to D$ by automorphisms.  

We take $p \in P$,
 and set $x :=\varpi(p) \in D$. 
If $h \in H_x$, 
 then $\varpi(h p) = h \cdot x = x = \varpi(p)$.
Therefore, 
there is a unique element of $K$,
 denoted by  $i_p(h)$,
 such that
\begin{equation}\label{eqn:3.1.3}
    h p = p \, i_p(h).  
\end{equation}
The correspondence $h \mapsto i_p(h)$
 gives rise to a Lie group homomorphism
$
    i_p:H_x \to K.
$
We set
\begin{equation}\label{eqn:3.1.4}
   H_{(p)}:=i_p(H_x).  
\end{equation}
Then, $H_{(p)}$ is a subgroup of $K$.

\begin{defin}
\label{def:autoP}
By an automorphism of the $H$-equivariant principal $K$-bundle
$\varpi : P \to D$,
we mean that there exist a diffeomorphism
 $\sigma: P \to P$ 
 and Lie group automorphisms $\sigma: K \to K$
and $\sigma: H \to H$ 
(by a little abuse of notation,
we use the same letter $\sigma$)
 such that
\begin{texteqn}
     $\sigma( h p k)= \sigma(h)\sigma(p)\sigma(k)
     \quad
     (h \in H, k \in K, p \in P).$
\label{eqn:3.1.1}
\end{texteqn}
\end{defin}
The condition \eqref{eqn:3.1.1} immediately implies:
\begin{texteqn}
     {$\sigma$ induces an action (denoted again by 
      $\sigma$) on $P/K \simeq D$,}
\label{eqn:sigmaD}
\end{texteqn}
\begin{texteqn}
     {the induced action $\sigma$ on $D$ is compatible with
      $\sigma \in \Aut(H)$ (see \eqref{eqn:compatHD} for the definition).}
\label{eqn:sigmaC}
\end{texteqn}

\medskip

We write $P^\sigma$ for the set of fixed points by $\sigma$, 
that is,
$$
P^\sigma := \set{p \in P}{\sigma(p) = p}.
$$

Then, we have: 
\begin{lem}
$\sigma(H_{(p)}) = H_{(p)}$
if $p \in P^\sigma$.
\end{lem}

\begin{proof}
Take $h \in H_x$.
Applying $\sigma$ to the equations $h \cdot x = x$ $(\in D)$
and $hp = pi_p(h)$ $(\in P)$,
we have $\sigma(h) \cdot x = x$
and $\sigma(h)p = p\sigma(i_p(h))$ 
from \eqref{eqn:3.1.1}.
Hence, $\sigma(h) \in H_x$
and $i_p(\sigma(h)) = \sigma(i_p(h))$.
Therefore, $\sigma(H_x) \subset H_x$ and
$\sigma(H_{(p)}) \subset H_{(p)}$.
Likewise,
$\sigma^{-1}(H_x) \subset H_x$ and
$\sigma^{-1}(H_{(p)}) \subset H_{(p)}$.
Hence, we have proved
$\sigma(H_x) = H_x$ and
$\sigma(H_{(p)}) = H_{(p)}$.
\end{proof}

\subsection{Multiplicity-free theorem}

For a representation $\mu$ of $K$,
 we denote by $\contra{\mu}$ the contragredient representation of $\mu$.
It is isomorphic to the conjugate representation $\overline{\mu}$
 if $\mu$ is unitary.
\begin{prop}\label{thm:3.2}  
Retain the setting of Subsection~\ref{sect:3.1}.  
Assume that there exist
 an automorphism $\sigma$ of 
the $H$-equivariant principal $K$-bundle
$\varpi : P \to D$
such that
\begin{align}
&\phantom{(a)}\parbox[t]{27em}{the induced action of $\sigma$ on $D$ 
is anti-holomorphic,}
\label{eqn:3.1.2}
\end{align}
and a subset $B$ of $P^\sigma$
 satisfying the following two conditions:
\begin{align}  
  &\phantom{(a)}\text{$H B K$ contains a non-empty open subset of $P$.  }
\label{eqn:3.2.1}
\\[\medskipamount]
  &\phantom{(a)}\parbox[t]{27em}{The restriction 
$\mu|_{H_{(b)}}$ is multiplicity-free as an $H_{(b)}$-module 
for any $b \in B$.}
\label{eqn:3.2.2}
\end{align}
We write its irreducible decomposition as
 $\mu|_{H_{(b)}}\simeq \bigoplus\limits_{i=1}^n \nu_b^{(i)}$. 
Further, we assume: 
\begin{subequations}\label{eqn:3.2.3}
  \renewcommand{\theequation}{\theparentequation) (\alph{equation}}%
\abovedisplayskip0pt
\begin{align}
 &\parbox[t]{27em}{${\mu\circ\sigma}\simeq \contra{\mu}$ as $K$-modules.}
\label{eqn:3.2.3a}
\\
  &\text{For any $b \in B$ and $i$,
    ${\nu^{(i)}\circ\sigma}\simeq \contra{\nu^{(i)}}$ as 
   $H_{(b)}$-modules.}
\label{eqn:3.2.3b}
\end{align}  
\end{subequations}
Then,
 any unitary representation of $H$ that is
 realized in $\OV$ is multiplicity-free.  
\end{prop}  
The proof of 
Proposition~\ref{thm:3.2}   
is given in Section~\ref{sect:5}

\begin{rem}  
Loosely,
 the conditions  
\eqref{eqn:3.2.1}  
and 
\eqref{eqn:3.2.2}  
mean that 
 the holomorphic bundle $\mathcal{V} \to D$
 cannot be `too large',
 with respect to the transformation group $H$.
The remaining condition \eqref{eqn:3.2.3}  
is
  often automatically fulfilled
(e.g.\ Corollary~\ref{crl:3.4}).
\end{rem}  

\begin{rem}\label{rem:3.2.2}  
As in Remark~\ref{rem:2.2.2},
 Proposition~\ref{thm:3.2}
still holds 
if $H_{(b)}$ is replaced by its arbitrary subgroup
 $H'_{(b)}$ for each $b \in B$ in 
\eqref{eqn:3.2.2}  
and 
\eqref{eqn:3.2.3b}.  
\end{rem}  

\begin{rem}
\label{rem:Weylinv}
For a connected compact Lie group $K$,
the condition \eqref{eqn:3.2.3a} is satisfied 
for any finite dimensional representation $\mu$ of $K$ 
if we take 
$\sigma \in \Aut(K)$ to be a \textit{Weyl involution}. 
We recall that $\sigma$ is a Weyl involution
 if there exists a Cartan subalgebra $\mathfrak{t}$ of the Lie
algebra $\mathfrak{k}$ of $K$ such that 
$d\sigma = -\operatorname{id}$ on $\mathfrak{t}$.
It is noteworthy that any simply-connected compact Lie group admits a
Weyl involution.
\end{rem}

\subsection{Multiplicity-free theorem (third form)}
In the assumption of Proposition~\ref{thm:3.2},  
 the subgroups $H_{(b)}$ may depend on $b$
(see \eqref{eqn:3.2.2} and \eqref{eqn:3.2.3b}).
For actual applications, 
we give a weaker but simpler form by taking
 just one subgroup $M$ instead of a family of
 subgroups $H_{(b)}$.

For a subset $B$ of $P$,
 we define the following subgroup $M_H(B)$ of $K$:
\begin{align}
\label{eqn:MHB}
   M_H(B)
:={}&\set{k \in K}{\text{ for each $b \in B$,
 there is $h \in H$ such that $h b = b k$}}  
\\
={}& \bigcap_{b \in B} K_{Hb} \; ,
\nonumber
\end{align}
where $K_{Hb}$ denotes the isotropy subgroup at 
$Kb$ in the left coset space $H\backslash P$, 
which is acted on by $K$ from the right.
Then $M_H (B)$ is $\sigma$-stable if $B \subset P^\sigma$,
 as is readily seen from \eqref{eqn:3.1.1}.

\begin{thm}\label{thm:3.3}  
Assume that there exist
 an automorphism $\sigma$ of the $H$-equivariant principal $K$-bundle
$\varpi: P \to D$
 satisfying 
\eqref{eqn:3.1.2}  
and a subset $B$ of $P^\sigma$
 with the following three conditions 
\eqref{eqn:3.3.1} -- \eqref{eqn:3.3.3}:
Let $M := M_H (B)$.
\begin{align}  
  &\phantom{(a)}\text{$H B K$ contains a non-empty open subset of $P$.  }
\label{eqn:3.3.1}
\\[\medskipamount]
 &\phantom{(a)}
\parbox[t]{27em}{The restriction $\mu|_M$ is multiplicity-free. }
\label{eqn:3.3.2}
\end{align}
 We shall write its irreducible decomposition as
 $\mu|_{M}\simeq \bigoplus\limits_{i=1}^n \nu^{(i)}$.
\begin{subequations}\label{eqn:3.3.3}
  \renewcommand{\theequation}{\theparentequation) (\alph{equation}}%
\abovedisplayskip0pt
\begin{align}
 &\parbox[t]{27em}
{${\mu\circ\sigma}\simeq \contra{\mu}$ as representations of $K$.}
\label{eqn:3.3.3a}
\\
  &\text{$\nu^{(i)}\circ\sigma \simeq \contra{\nu^{(i)}}$ as 
  representations of $M$
 for any $i$ $(1 \le i \le n)$.}
\label{eqn:3.3.3b}
\end{align}  
\end{subequations}
Then,
 any unitary representation 
 of $H$ which is realized in $\OV$ is multiplicity-free.  
\end{thm}  

\begin{rem}\label{rem:MHB}
Theorem~\ref{thm:3.3} still holds if we 
 replace $M$ with an arbitrary $\sigma$-stable subgroup of 
$M_H (B)$ to verify the conditions \eqref{eqn:3.3.2} and
 \eqref{eqn:3.3.3b}. 
\end{rem}

Assuming Proposition~\ref{thm:3.2},
we first complete the proof of Theorem~\ref{thm:3.3}.

\begin{proof}[Proof of Theorem~\ref{thm:3.3}]  
In view of Proposition~\ref{thm:3.2}   
and Remark~\ref{rem:3.2.2},  
it is sufficient to 
 show $M_H(B) \subset H_{(b)}$
 for all $b \in B$.

To see this,
take any $k \in M_H(B)$.
By the definition \eqref{eqn:MHB},
there exists $h \in H$ such that $hb = bk$.
Then, $h \in H_{\varpi(b)}$.
Since $i_b(h) \in K$ is characterized by the property 
$
    h b = b \, i_b(h)
$
(see \eqref{eqn:3.1.3}),
 $k$ coincides with $i_b(h)$. 
Hence,
$k = i_b(h) \in i_b(H_{\varpi(b)}) = H_{(b)}$
(see \eqref{eqn:3.1.4}).
Thus, we have proved
$M_H(B) \subset H_{(b)}$ for all $b \in B$.
\end{proof}  

\subsection{Line bundle case}
In general, 
the condition \eqref{eqn:3.3.1} tends to be fulfilled if $B$ is large, 
while the condition \eqref{eqn:3.3.2} tends to be fulfilled if $B$ is
small (namely, if $M$ is large).
However, we do not have to consider 
the condition \eqref{eqn:3.3.2} if 
 $\mathcal{V} \to D$ is a line bundle.
Hence, by taking $B$ to be maximal,
that is, by setting $B := P^\sigma$,
we get: 
\begin{crl}  
\label{crl:3.4}
Suppose we are in the setting of Subsection~\ref{sect:3.1}.
Suppose furthermore that $K$ is connected and $\dim \mu=1$.
Assume that there exists
 an automorphism $\sigma$ of the $H$-equivariant principal $K$-bundle
$\varpi: P \to D$
 satisfying 
\eqref{eqn:3.1.2}  
and the following two conditions:
\begin{align}  
  &\text{$d \sigma=-\operatorname{id}$
       on the center $\mathfrak{c}(\mathfrak{k})$ 
       of the Lie algebra $\mathfrak{k}$ of $K$.}
\label{eqn:3.4.1}
\\[\medskipamount]
  &\text{$H P^{\sigma}K$ contains a non-empty open subset of $P$.}
\end{align}  
Then, any unitary representation which can be realized in 
 $\OV$ is multiplicity-free.  
\end{crl}  
\begin{proof}[Proof of Corollary]  
As we mentioned, we apply Theorem~\ref{thm:3.3}   
with $B:=P^{\sigma}$.
The condition \eqref{eqn:3.3.2}   
is trivially satisfied
 because $\dim \mu =1$.  

Let us show $\mu \circ \sigma = \contra{\mu}$.
We write $K = [K,K]\cdot C$,
where $[K,K]$ is the commutator subgroup and
$C=\exp(\mathfrak{c}(\mathfrak{k}))$. 
Since $[K,K]$ is semisimple,
it acts trivially on the one dimensional representations
$\mu\circ\sigma$ and $\contra{\mu}$.
By \eqref{eqn:3.4.1},
$\mu\circ\sigma(e^X)=\mu(e^{-X})=\contra{\mu}(e^X)$
for any $X \in \mathfrak{c}(\mathfrak{k})$.
Hence $\mu\circ\sigma=\contra{\mu}$ both on $[K,K]$ and $C$.
Therefore, the condition \eqref{eqn:3.3.3a} holds.
Then, \eqref{eqn:3.3.3b}   
also holds.  
Therefore,
 Corollary follows from Theorem~\ref{thm:3.3}.     
\end{proof}  

\subsection{Multiplicity-free branching laws}
So far,
 we have not assumed that $P$ has a group structure.  
Now, we consider the case 
 that $P$ is a Lie group which we denote by $G$,
 and that $H$ and $K$ are closed subgroups of $G$.  
This framework enables us to apply Theorem~\ref{thm:3.3} to the
 \textit{restriction} of representations of $G$
(constructed on $G/K$)
to its subgroup $H$. 
Applications of Corollary~\ref{cor:4.6} include multiplicity-free
 branching theorems of highest weight representations 
for both finite and infinite dimensional cases 
(see \cite{xkglvis,xkrims40,xkmfkorea}).

We denote the centralizer of $B$ in $H \cap K$ by
$$
  Z_{H \cap K}(B):=\set{l \in H \cap K}{l b l^{-1}=
   b \text{ for any }b \in B}.
$$

\begin{crl}\label{cor:4.6}
Suppose $D = G/K$ carries a $G$-invariant complex structure,
and $\mathcal{V} = G \times_K V$
is a $G$-equivariant 
holomorphic vector bundle over $D$  associated to
a unitary representation $\mu: K \to GL(V)$. 
We assume there exist 
an automorphism $\sigma$ of the Lie group $G$ stabilizing $H$ and $K$
such that the induced action on $D = G/K$ is anti-holomorphic,
and
a subset $B$ of $G^\sigma$
satisfying the conditions \eqref{eqn:3.3.1}, \eqref{eqn:3.3.2}, and 
\eqref{eqn:3.3.3a} and {\upshape(b)} for
$P := G$ and $M := Z_{H \cap K} (B)$.
Then,
any unitary representation of $H$ which can be realized in 
the $G$-module
$\mathcal{O}(D,\mathcal{V})$
is multiplicity-free.
\end{crl}

\begin{proof}
Since $Z_{H\cap K}(B)$ is contained in $M_H(B)$ by the definition
\eqref{eqn:MHB}, 
Corollary~\ref{cor:4.6} is a direct consequence of
Theorem~\ref{thm:3.3} and Remark~\ref{rem:MHB}.
\end{proof}

\section{Proof of Proposition~\protect\ref{thm:3.2}}
\label{sect:5}

This section gives a proof of
Proposition~\ref{thm:3.2}    
by showing that
 all the conditions of Theorem~\ref{thm:mfvis} are
fulfilled. 
Then, the proof of our third form
(Theorem \ref{thm:3.3}) will be completed.

\subsection{Verification of the condition (\protect\ref{eqn:42B})}
Suppose we are in the setting of Proposition~\ref{thm:3.2}.
Then, 
 $HBK$ contains a non-empty open subset of $P$,
and consequently 
$\varpi(HBK)$ contains a non-empty open subset, say $W$, of $D$.
By taking the union of $H$-translates of $W$,
we get an $H$-invariant open subset
$D' := H \cdot W$ of $D$.
We set
$$
S := D' \cap \varpi(B).
$$
Then, $D' = H \cdot S$.
Besides,
$\sigma |_S = \operatorname{id}$
because $B \subset P^\sigma$.
Thus, the $H$-action on $D$ is $S$-visible with a compatible
automorphism $\sigma$ of $H$ by \eqref{eqn:3.1.1} 
in the sense of Definition~\ref{def:Scompat}.
Thus, the condition \eqref{eqn:42B} holds for $D'$.

\subsection{Verification of the condition (\protect\ref{eqn:42F})}
Next, let us prove
 that $\mathcal{V}_x$ is multiplicity-free as an
$H_x$-module for all $x$ in $S$.

Let $\mathcal{V}\simeq P \times_K V$ be the associated bundle, 
and
$
    P \times V \to \mathcal{V},
   \
   (p,v)\mapsto [p,v]
$
by the natural quotient map. 
For $p \in P$ we set $x :=\varpi(p) \in D$. 
Then, we can 
 identify the fiber $\mathcal{V}_x$
   with $V$ by the bijection 
\begin{equation}
\label{eqn:6.2.1}
   \iota_p:V \rarrowsim \mathcal{V}_x,
   \quad
   v \mapsto [p,v].  
\end{equation}

Via the bijection \eqref{eqn:6.2.1}
 and the group homomorphism $i_p :H_x \to H_{(p)}$,
 the isotropy representation of $H_x$ on $\mathcal{V}_x$
 factors through the representation $\mu: H_{(p)}\to GL(V)$,
 namely,
 the following diagram commutes for any $l \in H_x$:
\begin{equation}\label{eqn:3.7.2}
  \CD
        V \,\,@>{\sim}>{\iota_p}> \,\, \mathcal{V}_{x}
\\
        @V{\mu(i_p(l))}VV \,\,             \,\,@VV{L_l}V
\\
        V \,\,@>{\sim}>{\iota_p}>\,\,\mathcal{V}_x
  \endCD
\end{equation}

Now, suppose $x \in S$.
We take $b \in B$ such that $x = \varpi(b)$.

According to \eqref{eqn:3.2.2},  
 we decompose $V$ 
as a multiplicity-free sum of irreducible representations
 of $H_{(b)}$,
for which we write
\begin{equation}
\label{eqn:munub}
\mu = \bigoplus_{i=1}^n \nu_b^{(i)}, \quad  
  V = \bigoplus_{i=1}^n V_b^{(i)}.
\end{equation}
Then,
 it follows from \eqref{eqn:3.7.2}   
that if we set
$\mathcal{V}_{x}^{(i)}:=\iota_b (V_b^{(i)})$,
then
\begin{equation}\label{eqn:3.9.1}
     \mathcal{V}_{x}=\bigoplus_{i=1}^n \mathcal{V}_{x}^{(i)}
\end{equation}
is an irreducible decomposition
 as an $H_{x}$-module.  
Hence, \eqref{eqn:42F} is verified.

\subsection{Verification of the condition (\protect\ref{eqn:42C})}

Third, let us construct an isomorphism
$\Psi: \mathcal{V} \to \overline{\sigma^*\mathcal{V}}$.
According to the assumption \eqref{eqn:3.2.3a},  
 there exists a $K$-intertwining isomorphism,
 denoted by $\psi:V \to \overline V$, 
 between the two representations
 $\mu$ and $\overline{\mu \circ \sigma}$.
As the vector bundle $\mathcal{V} \to D$ is associated to the data 
$(P,K,\mu,V)$,
so is the vector bundle $\overline{\sigma^*\mathcal{V}} \to D$ to the data
$(P,K,\overline{\mu\circ\sigma}, \overline{V})$.
Hence  
 the map
$$
     P \times V \to P \times \overline V,
     \quad
     (p,v)\mapsto (p,\psi(v))
$$
 induces the bundle isomorphism
\begin{equation}
   \Psi:\mathcal{V} \rarrowsim \overline{\sigma^* \mathcal{V}}.
\end{equation}

In other words, 
the conjugate linear map defined by
\begin{equation}
\label{eqn:varphipsi}
\varphi: V \to V,
\quad v \mapsto \overline{\psi(v)}
\end{equation}
satisfies
\begin{equation*}
\mu(\sigma(k))\circ\varphi = \varphi\circ\mu(k)
\quad\text{for $k\in K$}.
\end{equation*}
Hence, we can define an anti-holomorphic endomorphism of $\mathcal{V}$
by 
\begin{equation*}
\mathcal{V} \to \mathcal{V}, 
\quad [p,v] \mapsto [\sigma(p), \varphi(v)].
\end{equation*}
This endomorphism, denoted by the same letter $\sigma$, 
 is a lift of the anti-holomorphic map
$\sigma: D \to D$, 
and satisfies \eqref{eqn:HisomV} because of \eqref{eqn:3.1.1}. 

Besides, for $x = \varpi(p)$,
we have
\begin{equation}
\label{eqn:iotasigma}
\iota_{\sigma(p)} \circ \varphi
= \sigma_x \circ \iota_p \; .
\end{equation}

Finally, let us verify the condition \eqref{eqn:42C}(a).
\par\noindent
\textbf{Step 1}.\enspace
 First, let us show 
\begin{equation}
\label{eqn:psitwo}
 \varphi(V_b^{(i)})={V_b^{(i)}}
  \quad \text{for $1 \le i \le n$. }
\end{equation}
Bearing the inclusion $H_{(b)} \subset K$ in mind, 
we consider the representation
$\overline{\mu\circ \sigma} : K \to GL(\overline{V})$ 
and its subrepresentation realized on 
$\psi(V_b^{(i)})$ $(\subset \overline{V})$
as an $H_{(b)}$-module.
Then, this is isomorphic to $(\nu_b^{(i)}, V_b^{(i)})$
 as $H_{(b)}$-modules because
 $\psi:V \to \overline V$
 intertwines the two representations $\mu$ and $\overline{\mu \circ \sigma}$
 of $K$. 
On the other hand,
it follows from the irreducible decomposition \eqref{eqn:munub} 
that the representation 
$\overline{\mu \circ \sigma}$
when restricted to the subspace
$\overline{V_b^{(i)}}$ is isomorphic to
$\overline{\nu_b^{(i)} \circ \sigma}$
as $H_{(b)}$-modules.
By our assumption \eqref{eqn:3.2.3b},  
$\nu_b^{(i)}$ is
 isomorphic to $\overline{\nu_b^{(i)}\circ \sigma}$,
 which occurs in $\overline V$
exactly once.
Therefore,
the two subspaces $\psi(V_b^{(i)})$ and $\overline{V_b^{(i)}}$
must coincide.  
Hence, we have \eqref{eqn:psitwo} by \eqref{eqn:varphipsi}.

\noindent
\textbf{Step 2}.\enspace
Next we show that \eqref{eqn:42C}(a)  
holds for $x = \varpi(b)$
if $b \in B$.  
We note that $\sigma(b)=b$ and $\sigma(x) = x$.
Then, 
it follows from \eqref{eqn:iotasigma}   
and \eqref{eqn:psitwo} 
that
\begin{equation*}
\sigma_x\circ\iota_b(V_b^{(i)})
= \iota_{\sigma(b)}\circ\varphi(V_b^{(i)})
= \iota_{\sigma(b)}(V_b^{(i)})
= \iota_b(V_b^{(i)}).
\end{equation*}
Since $\mathcal{V}_x^{(i)} = \iota_b(V_b^{(i)})$,
we have proved
$\sigma_x(\mathcal{V}_b^{(i)}) = \mathcal{V}_b^{(i)}$.

Hence, \eqref{eqn:42C}(a)  
holds. 

Thus, all the conditions of Theorem~\ref{thm:mfvis} hold for $D'$. 
Therefore, 
Proposition~\ref{thm:3.2} follows from Theorem~\ref{thm:mfvis}
and Remark~\ref{rem:2.2.2} (2).
Hence, the proof of Theorem \ref{thm:3.3} is completed.
\qed

\section{Concluding remarks}

\subsection{Applications in concrete settings}

The application of 
our multiplicity-free theorem ranges from finite dimensional
representations to infinite dimensional ones,
from discrete spectrum to continuous spectrum,
and from classical groups to exceptional groups.
Although concrete applications
 are not the main issue of the present article,
 let us mention some of them
 (see \cite{xkrims40, xkmfkorea} for details on this topic).  

The first paper \cite{xkmfjp} in this direction 
 (i.e. multiplicity-free theorem for the line bundle case) 
 already demonstrated that
there are fairly many \textit{new} multiplicity-free representations
for which explicit decomposition formulae were not known at that time.
(See \cite{Ali,Kra,Oka} and references therein in
the finite dimensional cases and to \cite{xkmfkorea,Seki} in the infinite
dimensional cases for some of new explicit branching laws.)

More generally, 
Theorem~\ref{thm:3.3}
 gives a systematic and synthetic proof of the multiplicity-free property 
including the Plancherel formula
for Riemannian symmetric spaces due to \'{E}. Cartan and
I. M. Gelfand,
its extension to line bundles and certain vector bundles
(A. Deitmar \cite{Deit}, see also \cite[Theorem 30]{xkrims40}),
and even its deformation which traces back to the
{\it{canonical representation}} of
Vershik--Gelfand--Graev in the $SL(2,\mathbb{R})$ case
(\cite{xdijk}, \cite[Example 8.3.3]{xkrims40} for more general groups);
the Hua--Kostant--Schmid $K$-type formula 
\cite{xkrims40,xschmidherm},
and its generalization to semisimple symmetric pairs due to the author
(\cite{xkmfjp}, see also \cite{xkmfkorea}).
Not only these infinite dimensional multiplicity-free
representations, 
but also quite a few multiplicity-free representations even in the
finite dimensional case where combinatorial argument is usually involved
 can be obtained as a special case
 of our propagation theorem.   
For example,
 we have proved in \cite{xkglvis}
 that Theorem \ref{thm:3.3} gives us 
a new and simple geometric explanation about when the tensor product of 
two representations for $GL(n)$ 
 is multiplicity-free (the classification was recently proved by 
Stembridge \cite{xstemgl} 
 in the spirit of case-by-case).

Least but not last,
 Theorem \ref{thm:3.3} has raised also a set of new problems concerning 
\textit{analysis on multiplicity-free representations}
 beyond (algebraic) branching laws, 
see \cite[Section 1.8]{xkmfkorea} for a short summary of developments
 made by Ben Sa\"id, 
 van Dijk, Hille, \O rsted, Neretin, 
 Zhang, 
 and by the author among others 
 in the last decade.

\subsection{Visible actions and coisotropic actions}

There are the following three concepts on group actions 
 in different geometric settings:
\begin{list}{}{\setlength{\itemsep}{0pt}
               \setlength{\leftmargin}{1.5em}}
\item[$\bullet$]
(Complex geometry) 
($S$-)visible actions (Definition \ref{def:S}).  
\item[$\bullet$]
(Symplectic geometry) 
coisotropic actions.
\item[$\bullet$]
(Riemannian geometry)
polar actions.
\end{list}
See \cite{xkglvis,xkrims40} for more details about visible actions
 on complex manifolds;
Guillemin and Sternberg
\cite{GS} or Huckleberry and Wurzbacher \cite{HW}
for coisotropic actions on symplectic manifolds; 
and Heintze et al \cite{xHei} or 
Podest\`{a}--Thorbergsson \cite{PT1}
for polar actions on Riemannian manifolds.
It should be noted that Lie groups $G$ are usually assumed to be compact
 for coisotropic actions
 and also for polar actions in the literature, 
 whereas we allow $G$ to be non-compact for visible actions in \cite{xkglvis,xkrims40}
 so that we can apply this concept to the study 
of infinite dimensional representations of $G$.  

We may compare the above three concepts assuming that the manifold is
 K\"{a}hler so that it is endowed with
complex, symplectic, and Riemannian structures simultaneously.
It should be noted that, according to \cite{HW}, J. Wolf first suggested
the terminology \lq\lq multiplicity-free actions" for 
coisotropic actions in the symplectic setting.
Further study on coisotropic actions and multiplicity-free representations
 of compact Lie groups may be found in \cite{HW}.  
The relation of visible actions with coisotropic actions and polar actions 
 is discussed in
\cite[Section~4]{xkrims40}.

A special case is given by linear actions on Hermitian vector
spaces. 
\begin{prop}
\label{prop:linear}
Suppose $\tau:G \to GL(V)$ is a unitary representation
 of a compact Lie group $G$
 on a finite dimensional Hermitian vector space $V$.  
Then the following two conditions {\rm{(i)}} and {\rm{(ii)}} are equivalent.  
Further,
 the condition {\rm{(iii)}} implies {\rm{(i)}} and {\rm{(ii)}}.    
\begin{enumerate}
\item[{\rm{(i)}}] $G$ acts strongly visibly on $V$
 as a complex manifold.  
\item[{\rm{(ii)}}]
$G$ acts coisotropically on $V$ 
 as a symplectic manifold.  
\item[{\rm{(iii)}}] $G$ acts polarly on $V$ as a Riemannian manifold.  
\end{enumerate}
\end{prop}
The equivalence (i) $\Longleftrightarrow$ (ii) follows from 
 the classification of multiplicity-free linear actions
by V. Kac (irreducible case), 
C. Benson--G. Ratcliff and A. Leahy (reducible cases; see \cite{BR} and
references therein),
Huckleberry and Wurzbacher \cite{HW},
 and the classification of multiplicity-free linear visible actions by A. Sasaki
\cite{Sa1,Sa3}.  
The last statement follows from Dadok \cite{Da}. 
The converse of the last statement does not hold.
A counterexample 
is the natural action of $U(3) \times Sp(n)$ on $\mathbb C^3 \otimes \mathbb C^{2n}$, which is not polar but is strongly visible and coisotropic.

\subsection{Generalization of the main theorem}

So far we have assumed that the base space $D$ and the fiber $V$
are finite dimensional,
 and discussed representations realized
 in the space of holomorphic sections
 for an equivariant holomorphic bundle ${\mathcal {V}} \to D$.  
It may be interesting to consider an analog of the propagation theorem
of multiplicity-free property (Theorem \ref{thm:2.2}) in a more general setting.  
Among others,
 we raise the following two cases for generalization.  
\vskip 1pc
1. Visible actions on infinite dimensional complex manifolds.

It is plausible that our framework 
 and its idea would work in the infinite dimensional settings
 by careful analysis
 (see \cite{Neeb}, for example). 
A generalization to infinite dimensional
 setting applies to the following objects:
\begin{enumerate}
\item[] the fiber $V$,
\item[] the complex manifold $D$ (base space), 
\item[] the group $G$.
\end{enumerate}
Here,
 we have in mind also an application to branching
problems  of representations of infinite dimensional Lie groups,
e.g., 
an infinite dimensional analogue of \cite[Theorem A]{xkmfkorea} for
those are constructed by a generalized Borel--Weil
theorem.

\vskip 1pc
2. Dolbeault cohomologies for equivariant holomorphic vector bundles.

The point here is to replace the space
$\mathcal{O}(D,\mathcal{V})$ of holomorphic sections by the Dolbeault
cohomology group $H^j(D,\mathcal{V})$.
In this setting,
 we highlight irreducible unitary representations
corresponding to a \lq\lq{geometric quantization}\rq\rq\ of elliptic orbits
$\mathcal{O}_\lambda$.
For real reductive groups,
these representations
 are realized in the Dolbeault cohomology groups for equivariant
holomorphic line bundles $\mathcal{L}_\lambda$ over
$\mathcal{O}_\lambda$,
giving the maximal globalization of Zuckerman derived functor modules
$A_{\mathfrak{q}}(\lambda)$.  
(The localization of their contragredient representations
 are  $K_\mathbb C$-equivariant sheaves of ${\mathcal {D}}$-modules
 supported on closed $K_{\mathbb{C}}$-orbits
 on the generalized flag variety $G_{\mathbb{C}}/Q$
 by the Hecht--Mili\v ci\'c--Schmid--Wolf duality theorem
 \cite{HMSW}.)
A generalization of our propagation theorem would yield an interesting family
of multiplicity-free branching problems,
see \cite[Conjecture 4.2]{Zuckerman60}.



\begin{thebibliography}{99}
%
\small
\parskip=0pt
\itemsep=0pt
%
\bibitem{Ali}
H. Alikawa,
Multiplicity-free branching rules for outer automorphisms of simple
Lie groups,
\textit{J. Math. Soc. Japan} \textbf{59} (2007), 151--177.

\bibitem{BR}
C. Benson and G. Ratcliff, On multiplicity free actions. Representations of 
real and $p$-adic groups, 221--304, Lect. Notes Ser. Inst. Math. Sci. Natl. 
Univ. Singap., 2, Singapore Univ. Press, Singapore, 2004.

\bibitem{Da}
J. Dadok, Polar coordinates induced by actions of compact Lie groups,
\textit{Trans. Amer. Math. Soc.} \textbf{288} (1985), 125--137.

\bibitem{Deit}
A. Deitmar,
Invariant operators on higher $K$-types.
\textit{J. Reine Angew. Math.} \textbf{412} (1990), 97--107.

\bibitem{xdijk}
G. van Dijk and S. C. Hille,
Canonical representations related to hyperbolic spaces,
\textit{J. Funct. Anal.}
\textbf{147}
(1997),
109--139.

\bibitem{xft}   
J. Faraut and E. G. F. Thomas, 
Invariant Hilbert spaces of holomorphic functions,
\textit{J. Lie Theory} \textbf{9} (1999),  383--402.

\bibitem{GS}
V. Guillemin and S. Sternberg,
Multiplicity-free spaces,
\textit{J. Differential Geom.}
\textbf{19}
(1984),
31--56.

\bibitem{HMSW}
H.Hecht, D. Mili\v ci\'c, W. Schmid, and J. A. Wolf,
Localization and standard modules for real semisimple Lie groups. I.
The duality theorem, 
\textit{Invent. Math. }
\textbf{90}
(1987), 
297--332.


\bibitem{xHei}
E. Heintze, R. S. Palais, and C.-L. Terng, G. Thorbergsson, Hyperpolar 
actions on symmetric spaces, in: Geometry, Topology, and Physics, Conf. 
Proc. Lecture Notes Geom. Topology, IV, International Press, Cambridge, MA, 
1995, pp. 214--245.

\bibitem{HW}
A. T. Huckleberry and T. Wurzbacher,
Multiplicity-free complex manifolds,
\textit{Math. Ann.}
\textbf{286}
(1990),
261--280.

\bibitem{xkobasho}
S. Kobayashi,
Irreducibility of certain unitary representations, 
\textit{J. Math. Soc. Japan}
\textbf{20} (1968), 638--642.

\bibitem{xkmfjp} 
T. Kobayashi,
Multiplicity-free theorem in branching problems of
unitary highest weight modules,
\textit{Proceedings of the Symposium on Representation Theory
held at Saga, Kyushu 1997} (ed. K. Mimachi), (1997),
\href{http://dml.ms.u-tokyo.ac.jp/PSRT/PSRT_18/PSRT_18_houkoku_009-017.pdf}{9--17}.

\bibitem{xkglvis}  
\bysame,  
Geometry of multiplicity-free representations of $GL(n)$,
visible actions on flag varieties, and triunity,
\textit{Acta Appl. Math.}
\textbf{81} (2004),
\href{http://dx.doi.org/10.1023/B:ACAP.0000024198.46928.0c}{129--146}.

\bibitem{xkrims40}  
\bysame,  
Multiplicity-free representations and visible actions on complex
manifolds, 
\textit{Publ. Res. Inst. Math. Sci.} \textbf{41} (2005), 
\href{http://dx.doi.org/10.2977/prims/1145475221}{497--549}
(a special issue of Publications of RIMS commemorating the fortieth
anniversary of the founding of the Research Institute for Mathematical
Sciences).

\bibitem{xkmfkorea} 
\bysame,  
Multiplicity-free theorems of the restriction of unitary highest
weight modules with respect to reductive symmetric pairs,
\textit{Progr. Math.} \textbf{255} Birkh\"{a}user (2007), 
\href{http://dx.doi.org/10.1007/978-0-8176-4646-2_3}{45--109}.

\bibitem{visiblesymm} 
\bysame,   
Visible actions on symmetric spaces,
\textit{Transform. Group} \textbf{12} (2007), 
\href{http://dx.doi.org/10.1007/s00031-007-0057-4}{671--694}.

\bibitem{xkgencar}  
\bysame,   
A generalized Cartan decomposition for the double coset space
$(U(n_1) \times U(n_2) \times U(n_3)) \backslash
U(n) / (U(p) \times U(q))$,
\textit{J. Math. Soc. Japan} \textbf{59} (2007), 
\href{http://dx.doi.org/10.2969/jmsj/05930669}{669--691}.

\bibitem{vermabranch}
\bysame,  
Restrictions of generalized Verma modules to symmetric pairs,
submitted, \href{http://arxiv.org/abs/1008.4544}{arXiv:1008.4544}.

\bibitem{Zuckerman60}
\bysame,  
Branching problems of Zuckerman derived functor modules,
to appear in Contemporary Mathematics, Amer. Math. Soc., 
In: Representation Theory and Mathematical Physics Conference
 in honor of Gregg Zuckerman's 60th birthday, 
 \href{http://arxiv.org/abs/1104.4399}{arXiv:1104:4399}.  

\bibitem{Kra}
C. Krattenthaler, 
Identities for classical group characters of nearly
rectangular shape, 
\textit{J. Algebra} 
\textbf{209}
 (1998),
 1--64.

\bibitem{xmackey}
G. W. Mackey,
Induced representations of locally compact groups I,
\textit{Annals of Math.}
\textbf{55} (1952),
101--139.

\bibitem{xneeb}
K.-H. Neeb,
On some classes of multiplicity free representations,
\textit{Manuscripta Math.}
\textbf{92} (1997),
389--407.

\bibitem{Neeb}
K-H. Neeb,
Towards a Lie theory of locally convex groups,
\textit{Jpn. J. Math.}
\textbf{1}
 (2006),  
291--468. 


\bibitem{Oka}
S. Okada,  
Applications of minor summation formulas to
rectangular-shaped representations of classical groups, 
\textit{J. Algebra}
\textbf{205}
(1998), 337--367. 

\bibitem{PT1}
F. Podest\`a and G. Thorbergsson,
Polar and coisotropic actions on K\"ahler manifolds,
\textit{Trans. Amer. Math. Soc.}
\textbf{354}
(2002),
1759--1781.

\bibitem{Sa1}
A. Sasaki,
Visible actions on irreducible multiplicity-free spaces,
\textit{Int. Math. Res. Not.} IMRN 2009, no. 18, 3445--3466.

\bibitem{Sa2}
\bysame,  
A characterization of non-tube type Hermitian symmetric spaces by
visible actions,
\textit{Geom. Dedicata} \textbf{145} (2010), 151--158.

\bibitem{Sa3}
\bysame,  
Visible actions on reducible multiplicitiy-free spaces,
\textit{Int. Math. Res. Not.} IMRN 2011, no. 4, 885--929.

\bibitem{xschmidherm}
W. Schmid,
Die Randwerte holomorpher Funktionen auf hermitesch symmetrischen
            R\"{a}umen,
\textit{Invent. Math.}
\textbf{9} (1969-70),
61--80.

\bibitem{Seki}
H. Sekiguchi,
Branching rules of Dolbealt cohomology groups over indefinite
Grassmannian manifolds,
\textit{Proc. Japan Acad.} Ser A. Math. Sci. \textbf{87} (2011), 31--34.

\bibitem{xstemgl}
J. R. Stembridge,  
Multiplicity-free products of Schur functions, 
\textit{Ann. Comb.} \textbf{5} (2001), 113--121.


\bibitem{xvu}
D. A. Vogan, Jr., Unitarizability of certain series of representations,
\textit{Ann. Math.} \textbf{120} (1984), 141--187.

\bibitem{xwal}
N. R. Wallach,  On the unitarizability of derived functor modules, 
\textit{Invent. Math.} \textbf{78} (1984),
  131--141.

\bibitem{xwolf}
J. Wolf,
Representations that remain irreducible on parabolic subgroups.
Differential geometrical methods in mathematical physics (Proc. Conf., 
Aix-en-Provence/Salamanca, 1979),
pp. 129--144,
\textit{Lecture Notes in Math.} \textbf{836},
Springer, Berlin, 1980.

\end{thebibliography}
\end{document}